arxiv-id: 0706.2638, pwd: tsy2d

\documentclass{article}
\usepackage{amsmath, amsfonts, amssymb, bbm}

\newtheorem{lemma}{Lemma}

\newtheorem{theorem}{Theorem}
\newtheorem{thmi}[theorem]{Theorem (Mellin Inversion)}
\newtheorem{thpf}[theorem]{Theorem (Plancherel Formula)}
\newtheorem{thcv}[theorem]{Theorem}
\newtheorem{thls}[theorem]{Theorem \cite{ls}}
\newtheorem{thasw}[theorem]{Theorem \cite{asw}}
\newcommand{\reff}[1]{(\ref{#1})}
\newcommand{\bC} {\mathbb{C}}

\newcommand{\bE} {\mathbb{E}}
\newcommand{\cF} {\mathcal{F}}
\newcommand{\cG} {\mathcal{G}}
\newcommand{\cH} {\mathcal{H}}

\newcommand{\cM} {\mathcal{M}}

\newcommand{\sP} {\mathsf{P}}

\newcommand{\bR} {\mathbb{R}}
\newcommand{\ve} {\varepsilon}
\newcommand{\vp} {\varphi}

\newcommand {\vk} {\varkappa}
\newcommand {\sgn} {\,\mbox{\scriptsize sgn}}
\begin{document}
\title{Some Applications of the Mellin Transform to Branching Processes}
\author{Wolfgang P. Angerer\\\\UAM Cuajimalpa\\Departamento de Matem\'aticas Aplicadas y Sistemas\\D.F. M\'exico}
\maketitle
\thispagestyle{empty}
\newpage
\thispagestyle{empty}
\begin{abstract}
\noindent We introduce a Mellin transform of functions which live on all of $\bR$ and discuss its applications to the limiting theory of Bellman-Harris processes, and specifically Luria-Delbr\"uck processes. More precisely, we calculate the life-time distribution of particles in a Bellman-Harris process from their first-generation offspring and limiting distributions, and prove a formula for the Laplace transform of the distribution of types in a Luria-Delbr\"uck process in the Mittag-Leffler limit. As a by-product, we show how to easily derive the (classical) Mellin transforms of certain stable probability distributions from their Fourier transform.\\

\noindent {\it AMS} (2000) {\it subject classification}. 44A05, 60J80.\\


\end{abstract}
\newpage
\section{Introduction}
In this note, we shall concern ourselves with some applications of the Mellin transform to the theory of branching processes; see Sections 4 and 5 for further motivation. We start with a definition of the Mellin transform for functions that live on all of $\bR$. Next, we give an overwiew of the classical definition of the Mellin transform and some of its properties. This is material which, except perhaps for the Plancherel formula, is well-known and might as well have fitted into an appendix. We prove our main theorem (which is actually nothing more than a formula to relate the Mellin transform of a function to its Laplace transform) in Section 3. The final two sections are applications of the main theorem.\\ 

To begin, let $X$ and $Y$ be two $\bR$-valued random variables such that $\sP(X \leq x) =: \cF(x)$ and $\sP(Y \leq y) =: \cG(y)$. Consider their product $XY$. Then, by total probability,
\begin{eqnarray}
\lefteqn{\sP(X Y \leq z) = \int_{-\infty}^{\infty} \sP(X Y \leq z \vert Y \in dy)\, \sP(Y \in dy)}\nonumber\\
&=&\int_{-\infty}^{0^-} \sP(X \geq z Y^{-1} \vert Y \in dy)\, d\cG(y)\nonumber\\
&& + \int_{0^+}^{\infty} \sP(X \leq z Y^{-1} \vert Y \in dy)\, d\cG(y) +  \sP(0 \leq z) \cdot g_0\:,\nonumber
\end{eqnarray}
if $\cG$ charges $\{0\}$. It follows that if $X$ and $Y$ are independent,
\[
\sP(X Y \leq z) = \int_{-\infty}^{0^-} \big(1 - \cF(z y^{-1})\big)\, dG(y) + \int_{0^+}^{\infty} \cF(z y^{-1})\, dG(y) + \sP(0 \leq z) \cdot g_0\:,
\]
and if they both have a density, then
\[
f \circledast g(z) := \int^{\infty}_0 \big[f(z y^{-1}) g(y) +  f(-z y^{-1}) g(-y) \big] \, \frac{dy}{y}
\]
is the density of the product $X Y$. $f \circledast g$ is called the Mellin convolution of $f$ and $g$. We would like to calculate its Mellin transform (see the following section for the definition and a number of properties of the Mellin transform). To that end, we define
\begin{equation}
(f \circledast g)^+(z) = f \circledast g(z) \cdot I_{\bR^+}(z)
\end{equation}
and
\begin{equation}
(f \circledast g)^-(z) = f \circledast g(-z) \cdot I_{\bR^-}(-z) = f \circledast g(-z) \cdot I_{\bR^+}(z)\:.
\end{equation}
Then both $(f \circledast g)^+$ and $(f \circledast g)^-$ live on $\bR^+$. Now we define the Mellin transform of $f \circledast g$ as the pair
\[
\cM[f \circledast g(x); s] := \big(\cM (f \circledast g)^+(s), \cM (f \circledast g)^-(s)\big)\:,
\]
of Mellin transforms of $f \circledast g^+$ and $f \circledast g^-$, and similarly for any other function that lives on $\bR$. (We reserve the symbol $\cM f$ for the Mellin transform of functions $f$ that live on $\bR^+$.) Now
\begin{eqnarray} 
\lefteqn{\cM (f \circledast g)^+(s) = \int_0^{\infty} (f \circledast g)^+ (z) \, z^{s - 1} \,dz = \int_0^{\infty} f \circledast g(z)\,z^{s - 1}\,dz}\nonumber\\
&=& \int_0^{\infty} \int_0^{\infty} \big[f(z y^{-1}) g(y) +  f(-z y^{-1}) g(-y) \big] \, \frac{dy}{y}\, z^{s - 1}\, dz \nonumber\\
&=& \int_0^{\infty} y^s g(y)\int_0^{\infty} f(z y^{-1}) \left(\frac{z}{y}\right)^{s - 1}\,d\left(\frac{z}{y}\right)\,\frac{dy}{y}\nonumber\\
&&+\: \int_0^{\infty} y^s g^-(y) \int_0^{\infty} f^-(z y^{-1})\,\left(\frac{z}{y}\right)^{s - 1}\,d\left(\frac{z}{y}\right)\,\frac{dy}{y}\:,\nonumber
\end{eqnarray}
so that
\begin{displaymath}
\cM(f \circledast g)^+(s) =  \cM f^+ (s)\,\cM g^+(s) + \cM f^-(s)\,\cM g^-(s) \:.\label{re}
\end{displaymath}
Similarly, one finds
\begin{displaymath}
\cM(f \circledast g)^-(s) = \cM f^+(s)\,\cM g^-(s) + \cM f^- (s)\,\cM g^+(s)\:.\label{im}
\end{displaymath}
Define a product $\odot$ of two Mellin transforms according to
\begin{eqnarray}
\lefteqn{\cM [f(x); s]\odot\cM [g(x); s]}\label{hypro}\\
&:=& \big(\cM f^+(s)\,\cM g^+(s) + \cM f^-(s)\,\cM g^-(s), \cM f^+(s)\,\cM g^-(s) + \cM f^-(s)\,\cM g^+(s)\big)\:.\nonumber
\end{eqnarray}
We have proved
\begin{thcv} \label{cv} Let $X$ and $Y$ be random variables with densities $f$ and $g$, respectively. Then
\begin{equation}
\cM[f \circledast g (x); s] = \cM[f(x);s] \odot \cM[g(x);s]\:.
\end{equation}
\end{thcv}
This is well-known for functions that live on $\bR^+$. The multiplication $\odot$ is the same as the multiplication rule for hyperbolic numbers \cite{sob}; it remains to see whether this observation has more than curiosity value.
\section{The Mellin Transform}
Let $f$ be a complex-valued function on $\bR^+ \cup \{0\}$ such that the integral
\[
\cM[f(x); s] := \int_0^{\infty} f(x) \, x^{s - 1}\, dx
\]
exists for all complex $s$ in the {\it fundamental strip} $\langle \alpha, \beta \rangle := \{s \in \bC : \alpha < \Re(s) < \beta \}$. Then $\cM[f(x); s]$ is called the {\it Mellin} transform of $f$ with respect to $s$. In what follows, we shall only be concerned with the Mellin transform of functions which live on $\bR^+$, and consequently write $\cM f$ or $\cM f(s)$ instead of $\cM[f(x); s]$ when this is feasible. The Mellin transform has a number of interesting properties. We start with
\begin{lemma}
\begin{equation}
\cM[f(\lambda x); s] = \lambda^{-s} \cM[f(x); s]\label{sc}
\end{equation}
and
\begin{equation}
\cM[f(x^{\mu}); s] = \mu^{-1} \cM[f(x); s/\mu]\:.
\end{equation}
\end{lemma}
{\it Proof}. Consider
\begin{eqnarray}
\lefteqn{\cM[f(\lambda x^{\mu}); s] = \int_0^{\infty} f(\lambda x^{\mu})\,x^{s - 1}\, dx}\nonumber\\
&=& \frac{1}{\mu} \int_0^{\infty} f(\lambda y)\,y^{s/\mu - 1/\mu}\, y^{1/\mu-1}\,dy = \frac{\lambda^{-s/\mu}}{\mu} \int_0^{\infty} f(x)\,x^{s/\mu - 1}\,dx\:.\nonumber
\end{eqnarray}
Now set $\lambda = 1$ or $\mu = 1$. \hfill $\Box$\\

This was easy. The following is a little harder:
\begin{thmi} Suppose $\gamma$ belongs to the fundamental strip of $f$. Then
\[
f(x) = \frac{1}{2\pi i}\int_{\gamma - i\infty}^{\gamma + i\infty} \cM f(s)\,x^{-s} \,ds \:.
\]
\end{thmi}
We first prove the following
\begin{lemma}\label{first}
\[
\frac{1}{2\pi i}\int_{\gamma - i\infty}^{\gamma + i\infty} \cM f(s)\,ds = f(1)\:.
\]
\end{lemma}
{\it Proof of the lemma}. Fix $\ve > 0$. Then, because $e^{s^2}$ is analytic in $s$ everywhere in the complex plane, we have
\[
\int_{\gamma - i\infty}^{\gamma + i\infty}  e^{\pi \ve^2 s^2}\,ds = \int_{ - i\infty}^{i\infty}  e^{\pi \ve^2 s^2} \, ds = i \int_{ - \infty}^{\infty}  e^{- \pi \ve^2 s^2} \, ds = \frac{i}{\ve}\:,
\]
and so
\begin{eqnarray}
\lefteqn{\frac{1}{2\pi i}\int_{\gamma - i\infty}^{\gamma + i\infty}  e^{\pi \ve^2 s^2} \cM f(s)\,ds = \frac{1}{2\pi i}\int_{\gamma-i\infty}^{\gamma + i\infty}  e^{\pi \ve^2 s^2} \int_0^{\infty} f(z) \, z^{s - 1}\,dz\,ds}\nonumber\\
&=& \frac{1}{2\pi i}\int_{\gamma - i\infty}^{\gamma + i\infty}  e^{\pi \ve^2 s^2} \int_{-\infty}^{\infty} f(e^{-y}) \, e^{-ys} \, dy \,ds\nonumber\\
&=& \frac{1}{2\pi i} \int_{- \infty}^{\infty} f(e^{-y}) \int_{\gamma - i\infty}^{\gamma + i\infty}  e^{\pi \ve^2 s^2} e^{-2\pi \ve s y (2\pi \ve)^{-1}}\,ds\,dy\nonumber\\
&=& \frac{1}{2\pi \ve} \int_{-\infty}^{\infty} f(e^{-y}) \, e^{-\pi y^2  (2 \pi \ve)^{-2}}\, dy = \int_{-\infty}^{\infty} f(e^{-2 \pi \ve y}) \, e^{-\pi y^2}\, dy\:.\nonumber
\end{eqnarray}
Now let $\ve \to 0$.\hfill $\Box$\\

\noindent {\it Proof of the theorem}. By the lemma and the scaling property \reff{sc}, we have
\[
\frac{1}{2\pi i}\int_{\gamma - i\infty}^{\gamma + i\infty} \cM f(s)\,x^{-s} \,ds = \frac{1}{2\pi i}\int_{\gamma - i\infty}^{\gamma + i\infty} \cM [f(xy); s] \,ds = f(xy)\big\vert_{y = 1} = f(x)\:,
\]
which concludes the proof.\hfill $\Box$\\

The following theorem is in a similar spirit.
\begin{thpf} Suppose $\gamma$ belongs to the fundamental strip of both $f$ and $g$. Then
\[
\frac{1}{2\pi i}\int_{\gamma - i\infty}^{\gamma + i\infty} \cM f(s) \, \overline{\cM g(s)} \,ds = \int_{0}^{\infty} f(x)\,\overline{g(x)}\,x^{2 \gamma - 1}\,dx\:.
\]
\end{thpf}
{\it Proof}. As for the proof of Lemma \reff{first}, one checks that
\[
\frac{1}{2\pi i}\int_{\gamma - i\infty}^{\gamma + i\infty} \overline{\cM f(s)}\,ds =  \overline{f(1)}
\]
and
\[
\frac{1}{2\pi i}\int_{\gamma - i\infty}^{\gamma + i\infty} \overline{\cM f(s)} \, x^{-s} \,ds = \frac{x^{- 2 \gamma}}{2\pi i}\int_{\gamma - i\infty}^{\gamma + i\infty} \overline{\cM f(s)} \left(\frac{1}{x}\right)^{-\overline{s}} \,ds = x^{- 2 \gamma}\,\overline{f(1/x)}\:.
\]
Then
\begin{eqnarray}
\lefteqn{\frac{1}{2\pi i}\int_{\gamma - i\infty}^{\gamma + i\infty}  e^{\pi \ve^2 s^2} \cM f(s) \, \overline{\cM g(s)} \,ds}\nonumber\\
&=& \frac{1}{2\pi i}\int_{\gamma-i\infty}^{\gamma + i\infty}  e^{\pi \ve^2 s^2} \int_0^{\infty} f(x)\, x^{s - 1}\,dx \int_0^{\infty} y^{-2 \gamma} \overline{g(1/y)} \,y^{s - 1}\,dy\,ds\nonumber\\
&=& \frac{1}{2\pi i}\int_{\gamma-i\infty}^{\gamma + i\infty}  e^{\pi \ve^2 s^2} \int_{-\infty}^{\infty} f(e^{-u})\, e^{-us}\,du \int_{-\infty}^{\infty} e^{-2 \gamma v} \overline{g(e^{-v})} \,e^{vs}\,dv\,ds\nonumber\\
&=& \frac{1}{2\pi i} \int_{- \infty}^{\infty} f(e^{-u})  \int_{- \infty}^{\infty} e^{-2 \gamma v} \,\overline{g(e^{-v})} \int_{\gamma - i\infty}^{\gamma + i\infty}  e^{\pi \ve^2 s^2} e^{-2 \pi \ve s(u - v)(2 \pi \ve)^{-1}}\,ds\,du\,dv\nonumber\\
&=& \frac{1}{2\pi \ve} \int_{- \infty}^{\infty} e^{-2 \gamma v} \,\overline{g(e^{-v})} \int_{- \infty}^{\infty} f(e^{-(u - v) - v})\,e^{-\pi(u - v)^2(2 \pi \ve)^{-2}}\,du\,dv\nonumber\\
&=& \int_{- \infty}^{\infty} e^{-2 \gamma v} \,\overline{g(e^{-v})} \int_{- \infty}^{\infty} f(e^{-2 \pi \ve u - v})\,e^{-\pi u^2}\,du\,dv\:.\nonumber
\end{eqnarray}
Now let $\ve \to 0$ and change variables.\hfill $\Box$\\

We now continue with the derivation of
\section{The Main Theorem}
The gist of the paper is the following
\begin{theorem} \label{main} Let $f$ be a function on $\bR$ with bilateral Laplace transform
\[
\vp(u) := \int_{-\infty}^{\infty} e^{-ux} f(x)\,dx
\]
and Mellin transform $\cM[f(x); s] = \big(\cM f^+(s), \cM f^-(s)\big)$. Define
\begin{equation}
\alpha^+ := \sup \: \{\alpha \in \bR : \lim_{x \to \infty} e^{\alpha x} f^+(x) = 0\}
\end{equation}
and
\begin{equation}
\alpha^- := \sup \: \{\beta \in \bR : \lim_{x \to \infty} e^{\alpha x} f^-(x) = 0\}\:.
\end{equation}
Suppose that $\alpha^+$ and $\alpha^-$ are strictly larger than zero, and choose for $\beta$ any number in $\bR^-$ such that $0 < - \beta < \min (\alpha^+,\alpha^-)$. Then we have, for $s > 0$ and in the fundamental strip of $f^+$ and $f^-$, 
\begin{equation}
\cM f^+(s) = \frac{\Gamma(s)}{2 \pi i}\int_{\beta - i\infty}^{\beta + i\infty} \vp(u)(-u)^{-s}\,du \label{+m}
\end{equation}
and
\begin{equation}
\cM f^-(s) = \frac{\Gamma(s)}{2 \pi i}\int_{\beta - i\infty}^{\beta + i\infty} \vp(-u)(-u)^{-s}\,du\:.\label{-m}
\end{equation}
Furthermore,
\begin{eqnarray}
\lefteqn{\vp(u) = \frac{1}{2 \pi i}\int_{\gamma - i\infty}^{\gamma + i\infty} \cM f^+(s)\,\Gamma(1 - s)\,u^{s - 1}\,ds}\nonumber\\
&+&  \frac{1}{2 \pi i}\int_{\gamma - i\infty}^{\gamma + i\infty} \cM f^-(s)\,\Gamma(1 - s)\,(-u)^{s - 1}\,ds\label{inv}
\end{eqnarray}
for any $\gamma$ that belongs to the fundamental strips of both $f^+$ and $f^-$.
\end{theorem}
{\it Proof}. We first prove the theorem under the assumption that $f$ is concentrated on $\bR^+$, such that $f = f^+$, and we may write $\cM f$ for its Mellin transform. Choose $\beta := \beta^+$ as required. By Laplace inversion, we have
\[
f(x) = \frac{1}{2 \pi i}\int_{\beta - i\infty}^{\beta + i\infty} \vp(u) \, e^{u x}\,du
\]
and so
\[
2 \pi i \cM f(s) = \int_0^{\infty} \int_{\beta - i\infty}^{\beta + i\infty} \vp(u)\,e^{u x}\,du \,x^{s - 1} \, dx\:.
\]
The standard trick here is often to first convolve $f$ with some $(\pi \ve)^{-1/2}\,e^{-x^2/\ve}$. The Laplace transform of the latter is $e^{\ve u^2/4}$, and
\begin{eqnarray}
\lefteqn{\int_0^{\infty} \int_{\beta - i\infty}^{\beta + i\infty} e^{\ve u^2/4}\,\vp(u)\,e^{u x}\,du \,x^{s - 1} \, dx}\nonumber\\
&=& \int_{-\beta - i\infty}^{-\beta + i\infty} e^{\ve u^2/4}\,\frac{\vp(- u)}{u^s} \int_0^{\infty} e^{- x}\,x^{s - 1} \, dx\,du\nonumber\\
&=& \Gamma(s) \int_{\beta - i\infty}^{\beta + i\infty} e^{\ve u^2/4}\,\vp(u)\,(-u)^s \,du \nonumber\:,
\end{eqnarray}
where interchanging the order of integration is permissible because the integrand becomes exponentially small if either $u$ or $x$ are of the order $\sim \ve^{-1}$. One can then undo the convolution by letting $\ve \to 0$. As for the converse, choose $\gamma$ in the fundamental strip of $f$. Then
\[
f(x) = \frac{1}{2 \pi i}\int_{\gamma - i\infty}^{\gamma + i\infty}\cM f(s)\,x^{-s}\,ds\:,
\]
and
\begin{eqnarray}
\lefteqn{\int_0^{\infty}\int_{\gamma - i\infty}^{\gamma + i\infty} \cM f(s)\,x^{-s}\,ds\,e^{-ux}\,dx }\label{fc}\\
&=& \int_{\gamma - i\infty}^{\gamma + i\infty} \cM f(s)\, u^{s - 1}\int_0^{\infty} x^{-s}\,e^{-x}\,dx\,ds = \int_{\gamma - i\infty}^{\gamma + i\infty} \cM f(s)\,\Gamma(1 - s)\,u^{s - 1}\,ds \nonumber\:.
\end{eqnarray}
To justify the interchange of the order of integration, one can Mellin convolve $f$ with $(\pi \ve)^{-1/2}\,e^{-(\log x)^2/\ve}$, whose Mellin transform is $e^{\ve s^2/4}$, and then again undo the convolution by letting $\ve \to 0$. This proves the theorem in case $f$ is concentrated on $\bR^+$. For the general case, let $\vp^+$ and $\vp^-$ be the Laplace transforms of $f^+$ and $f^-$, respectively. Then
\begin{equation}
\vp(u) = \int_0^{\infty} f^+(x)\,e^{-u x}\,dx + \int_0^{\infty} f^-(x)\,e^{u x}\,dx = \vp^+(u) + \vp^-(-u)\:,\label{secc}
\end{equation}
and
\[
\cM f^+(s) = \frac{\Gamma(s)}{2 \pi i}\int_{\beta - i\infty}^{\beta + i\infty} \vp(u)\,(-u)^{-s}\,du - \frac{\Gamma(s)}{2 \pi i}\int_{\beta - i\infty}^{\beta + i\infty} \vp^-(-u)\,(-u)^{-s}\,du\:,
\]
by what we have proved already. Now $\vp^-$ cannot have any singularities in the half-plane $\Re(s) > - \beta$. We therefore can close the contour of integration in the second integral along a semi-circle of radius $r$ and find, by analyticity of $\vp^-$,
\[
\int_{\beta - i r}^{\beta + i r} \vp^-(-u)\,(-u)^{-s}\,du = - \: r \int_{-\pi/2}^{\pi/2} \vp^-(r e^{i \theta})\,(r e^{i \theta})^{-s}\,d\theta\:,
\]
which is of order
\[
r^{1 - s} \int_{-\pi/2}^{\pi/2} e^{-r \cos \theta}\,d\theta \leq 2 r^{1 - s} \int_0^{\pi/2} e^{-r (1 - 2\theta/\pi)}\,d\theta =  \pi r^{- s} (1 - e^{-r})
\]
in absolute value, and therefore tends to 0 if $s > 0$. This proves \reff{+m}. \reff{-m} is proved in a similar manner, and \reff{inv} is an easy consequence of \reff{fc} and \reff{secc}.\hfill$\Box$\\

Although we shall not use it as heavily, we note the following variant of Theorem \ref{main}:
\begin{theorem} Let $f$ be a function on $\bR$ with Mellin transform $(\cM f^+, \cM f^-)$ and Fourier transform
\begin{equation}
f^{\ast}(y) := \int_{-\infty}^{\infty} e^{i y x} f(x)\,dx =  \int_0^{\infty} e^{-i y x} f^-(x)\,dx + \int_0^{\infty} e^{i y x} f^+(x)\,dx\:.
\end{equation}
Then, for $0 < s < 1$,
\begin{equation}
\cM f^+ (s) = \frac{\Gamma(s)}{\pi}\,\Re \left(\exp\left(- i \, \frac{s\pi}{2} \right)\int_0^{\infty} f^{\ast}(y)\, y^{-s}\,dy\right)\label{fcs}
\end{equation}
and
\begin{equation}
\cM f^- (s) = \frac{\Gamma(s)}{\pi}\,\Re \left(\exp\left(i \, \frac{s\pi}{2} \right)\int_0^{\infty} f^{\ast}(y)\, y^{-s}\,dy\right)\:.\label{scs}
\end{equation}
Moreover,
\begin{eqnarray}
\lefteqn{f^{\ast}(y) = i\int_{\gamma - i\infty}^{\gamma + i\infty} \cM f^+(s)\,\Gamma(1 - s) \, \exp \left(- \frac{s \pi i}{2} \right)\,y^{s - 1}\,ds}\nonumber\\
&-&i \int_{\gamma - i\infty}^{\gamma + i\infty} \cM f^- (s)\,\Gamma(1 - s) \, \exp \left(\frac{s \pi i}{2} \right)\,y^{s - 1}\,ds \label{tcs}
\end{eqnarray}
for any $\gamma$ that belongs to the fundamental strips of $f^+$ and $f^-$.
\end{theorem}
{\it Proof}. We need the following
\begin{lemma}\label{cs} For $0 < s < 1$,
\begin{equation}
\int_0^{\infty} e^{\pm i x} \, x^{s-1}\,dx = \Gamma(s) \exp \left( \pm \frac{s \pi i}{2} \right)\:.
\end{equation}
\end{lemma}
{\it Proof of the lemma}. We have
\[
\int_0^{\infty} e^{\pm i x} e^{- x u}\,dx = \frac{1}{u \mp i} = \frac{u}{u^2 + 1} \pm \frac{i}{u^2 + 1}\:, 
\]
and then, by Theorem \ref{main},
\[
\int_0^{\infty} \cos x \,e^{- \lambda x} x^{s-1}\,dx = \frac{\Gamma(s)}{2 \pi i}\int_{\gamma - i\infty}^{\gamma + i \infty} \frac{u + \lambda}{(u + \lambda)^2 + 1} (-u)^s \, du 
\]
for arbitrary $\lambda > 0$. It then follows by the Residue theorem that
\[
\int_0^{\infty} \cos x \,e^{- \lambda x} \,x^{s-1}\,dx = \Gamma(s) \big( (\lambda - i)^{-s} + (\lambda + i)^{-s}\big) = \Gamma(s) \frac{\cos\big(- s \arg(\lambda - i)\big)}{(\lambda^2 + 1)^{s/2}}\:,
\]
which implies
\[
\int_0^{\infty} \cos x \, x^{s-1}\,dx = \Gamma(s) \cos \frac{s \pi}{2}\:,
\]
by arbitrariness of $\lambda$. Similarly, one proves
\[
\int_0^{\infty} \sin x \, x^{s-1}\,dx = \Gamma(s) \sin \frac{s \pi}{2}\:,
\]
and this already implies the lemma.\hfill $\Box$\\

\noindent {\it Proof of the theorem}. Inverting the Fourier transform, we find
\[
2 \pi f(x) = \int_{-\infty}^{\infty} f^{\ast}(y) \, e^{-i x y}\, dy\:,
\]
and so
\begin{eqnarray}
\lefteqn {2 \pi \cM f^+(s) = \int_0^{\infty} \int_{-\infty}^{\infty} f^{\ast}(y) \, e^{-i x y}\,dy \,x^{s - 1}\,dx}\nonumber\\
&=& \int_0^{\infty} f^{\ast}(-y) \int^{\infty}_0 e^{i x y}\,x^{s - 1}\,dx\,dy + \int_0^{\infty} f^{\ast}(y) \int_0^{\infty} e^{-i x y}\,x^{s - 1}\,dx\,dy\nonumber\\
&=& \Gamma(s)\exp\left(i\,\frac{s\pi}{2} \right)\int_0^{\infty} f^{\ast}(-y)\,y^{-s}\,dy + \Gamma(s)\exp\left(- i\,\frac{s\pi}{2} \right)\int_0^{\infty} f^{\ast}(y)\,y^{-s}\,dy\:,\nonumber
\end{eqnarray}
by Lemma \ref{cs}. But this already implies the first of Equations \reff{fcs} - \reff{tcs}, because $f^{\ast}(-y) = \overline{f^{\ast}(-y)}$. (One again checks that interchanging the order of integration is permissible at worst after convoluting $f$ with some $(\pi \ve)^{-1/2}\,e^{-x^2/\ve}$.) The rest are proved in a similar manner.\hfill$\Box$\\

As an application, we now calculate the Mellin transforms of those stable probability distributions with Fourier transform
\[
\psi^{\ast}(y) :=: \psi^{\ast}_{\alpha, \theta}(y) := \exp \left(- \vert y\vert^{\alpha} e^{i \pi \theta \sgn y /2} \right)\:.
\]
Writing $e^{i \pi \theta /2} =: \zeta$, we find
\begin{eqnarray}
\lefteqn{\int_0^{\infty} \psi^{\ast}(y)\,y^{-s}\,dy = \int_0^{\infty} \exp \left(- y^{\alpha} \zeta \right)\,y^{-s}\,dy}\nonumber\\
&=& \frac{1}{\alpha \zeta^{(1 - s)/\alpha}} \int_0^{\infty} e^{-y}\,y^{(1 - s)/\alpha - 1}\,dy = \Gamma\big((1 - s)/\alpha\big) \frac{\overline{\zeta}^{\,(1 - s)/\alpha}}{\alpha \vert \zeta\vert^{(1 - s)/\alpha}}\nonumber\\
&=& \alpha^{-1}\,\Gamma\big((1 - s)/\alpha\big) \left(\cos \frac{(1 - s)\pi \theta}{2\,\alpha} + i \sin \frac{(1 - s)\pi \theta}{2\,\alpha} \right)\:,\nonumber
\end{eqnarray}
because $\vert \zeta \vert = 1$. Multiplication by $e^{- i s\pi/2}$ and taking the real part then gives
\begin{eqnarray}
\lefteqn{\cM\psi^+(s) = \frac{\Gamma(s)\Gamma\big((1 - s)/\alpha\big)}{\alpha\,\pi}\cos \frac{(1 - s)\pi \theta + s \pi \alpha}{2\,\alpha}}\label{19}\\
&=& \frac{\Gamma(s)\Gamma\big((1 - s)/\alpha\big)}{\alpha\,\pi} \sin \frac{(1 - s)\pi (\alpha - \theta)}{2\,\alpha} = \frac{\rho^+\, \Gamma(s)\Gamma\big(1 + (1 - s)/\alpha\big)}{\Gamma\big(1 + \rho^+ (1 - s)\big)\Gamma\big(1 - \rho^+ (1 - s)\big)}\:,\nonumber
\end{eqnarray}
if we define
\begin{equation}
\rho^+ = \frac{\alpha - \theta}{2\,\alpha}\label{rho+}\:,
\end{equation}
and make use of $\Gamma(s)\,\Gamma(1 - s) = \pi (\sin \pi s)^{-1}$. \reff{19} is equivalent to Formula 17 in \cite{mpg} or Formula 6.8 in \cite{mlp}. Similarly, one finds
\begin{equation}
\cM \psi^-(s) = \frac{\rho^-\, \Gamma(s)\Gamma\big(1 + (1 - s)/\alpha\big)}{\Gamma\big(1 + \rho^- (1 - s)\big)\Gamma\big(1 - \rho^- (1 - s)\big)}\:,
\end{equation}
if we define
\begin{equation}
\rho^- = \frac{\alpha + \theta}{2\,\alpha}\label{rho-}\:.
\end{equation}
In particular, $\cM \psi^+(1) + \cM \psi^-(1) = \rho^+ + \rho^- = 1$, as required.
\section{An Application to Bellman-Harris Processes}
Let $\{Z_t\}_{t\geq 0}$ be a supercritical Bellman-Harris process with offspring distribution $\{\pi_k\}_{k = 0}^{\infty}$ and life-time distribution $G$. Thus, $G(t)$ is the probability that a newborn individual survives at least until time $t$, and $\pi_k$ is the probability that once it splits into a number $Z_+$ of progeny, it will split into exactly $k$ of these. We denote by $f(s) := E(e^{sZ_+}) = \sum_{k = 0}^{\infty} \pi_k s^k$ the corresponding generating function. By supercriticality, $\mu := f'(1) > 1$, and there exists $q \in [0, 1)$ such that $f(q) = q$. We assume that $G$ is non-lattice, that $G(0^+) = G(0) = 0$, and that $\mu > \infty$. Then there exist \lq constants' $\chi_t$ (the Seneta constants) such that, on the set of non-extinction, $\chi_t Z_t$ converges almost surely to a non-degenerate random variable $Z$ whose Laplace transform $\psi(u) = \bE(e^{-uZ})$ satisfies
\begin{equation}
\psi(u) = \int_0^{\infty} f \circ \psi(u e^{-\beta t})\,dG(t)\:,\label{leq}
\end{equation}
where $\beta$ is the Malthusian parameter, that is, the unique root of
\[
\mu \int_0^{\infty} e^{-y t}\,dG(t) = 1\:,
\]
$y \in (0, \infty)$. There are only a few instances where the solution of \reff{leq} is known for given $f$ and $G$. We turn the problem on its head: Suppose $\psi$ were (the Laplace transform of some function $g$) such that it fulfills the conditions of Theorem \ref{main}. Then
\[
\cM g(s)  = \frac{\Gamma(s)}{2 \pi i}\int_{\beta-i \infty}^{\beta + i\infty} \psi(u)\,(-u)^{-s}\,du
\]
for some suitably chosen value of $\beta$, and it follows by Equation \reff{leq} that
\begin{eqnarray}
\lefteqn{\frac{2\pi i}{\Gamma(s)}\,\cM g(s) = \int_{\beta-i \infty}^{\beta + i\infty} \int_0^{\infty} f \circ \psi(u e^{-t})\,dG(t) \,(-u)^{-s}\,du}\nonumber\\
&=& \int_0^{\infty} \int_{\beta-i \infty}^{\beta + i\infty} f \circ \psi(u e^{-t})\,(-u e^{-t})^{-s} e^{-t s}\,e^{t} d(ue^{-t})\,dG(t)\nonumber\\
&=& \int_{\beta-i \infty}^{\beta + i\infty} f \circ \psi(u)\,(-u)^{-s}\,du \int_0^{\infty} e^{-t (s-1)}\,dG(t)\:.\nonumber
\end{eqnarray}
In other words, we have for the Laplace transform of $G$ at $s$,
\begin{displaymath}
\int_0^{\infty} e^{-t s}\,dG(t) = \frac{\displaystyle \int_{\beta-i \infty}^{\beta+ i\infty} \psi(u)\,(-u)^{-s-1}\,du}{\displaystyle \int_{\beta-i \infty}^{\beta + i\infty} f \circ \psi(u)\,(-u)^{-s-1}\,du}\:.
\end{displaymath}
{\it Example}. Suppose that $g(x) \propto e^{-x} x^{\kappa - 1}$ for some $\kappa > 0$ and $f(s) = s^m$. Then
\begin{eqnarray}
\lefteqn{\int_{\beta-i \infty}^{\beta+ i\infty} \psi(u)\,(-u)^{-s-1}\,du = 2 \pi i \,\frac{\cM g(s + 1)}{\Gamma(s + 1)}}\nonumber\\
&=& \frac{2 \pi i }{\Gamma(s + 1)\,\Gamma(\kappa)}\,\int_0^{\infty} e^{-x}\,x^{\kappa - 1}\,x^s\,dx = 2 \pi i \,\frac{\Gamma(s + \kappa)}{\Gamma(s + 1)\,\Gamma(\kappa)}\:,\nonumber
\end{eqnarray}
which eventually implies that
\[
\int_0^{\infty} e^{-t s}\,dG(t) = \frac{\Gamma(s + \kappa)\,\Gamma(m\kappa)}{\Gamma(s + m\kappa)\,\Gamma(\kappa)}\:,
\]
and that $G$ has a density:
\begin{equation}
\frac{dG(t)}{dt} =  \frac{\Gamma(m\kappa)}{\Gamma(\kappa)\,\Gamma(m\kappa - \kappa)} \, e^{-\kappa t} (1 - e^{-t})^{(m - 1)\kappa - 1}\:.\label{pm}
\end{equation}
In case $f$ is a polynomial $f(s) = \sum_{j = 1}^m \pi_j s^j$, we similarly obtain
\begin{equation}
\int_0^{\infty} e^{-t s}\,dG(t) = \frac{\displaystyle \frac{\Gamma(s + \kappa)}{\Gamma(\kappa)}}{\displaystyle \sum_{i = 1}^m \pi_j \frac{\Gamma(s + j\kappa)}{\Gamma(j\kappa)}}\:.\label{poly}
\end{equation}
It is not obvious that the function on the right-hand side of \reff{poly} has all the properties of a Laplace transform (most of all, complete monotonicity), but this is what we have shown. So the question is, When does the Laplace transform $\psi$ of the random variable $Z$ have an exponential tail? The answer is given by the following
\begin{theorem} \label{expmo} Let $F_t$ be the PGF of particle numbers in a Bellman-Harris process at time $t$, and let $f$ be the PGF of the corresponding first-generation offspring distribution. Say that $f$ has exponential moments up to order $r > 0$ if $f(e^u) < \infty$ for $u < r$, and let $M_t := F'_t(1)$. Then $F_t$ has exponential moments up to order $O(r/M_t)$. In particular, there exists $r' > 0$ such that Laplace transform
\[
\tilde{\psi}(u) := \lim_{t \to \infty} F_t(e^{-u/M_t})
\]
is analytic for $u > -r'$.
\end{theorem}
{\it Proof}. We can assume $f(0) = 0$, so $Z_t \to \infty$ almost surely. Because
\[
e^u = F_t \circ \exp \left(\frac{M_t \log F_{-t}(e^u)}{M_t}\right)
\]
(we write $F_{-t}$ to denote the inverse of $F_t$), one readily checks that Theorem \ref{expmo} holds true iff
\[
\frac{d}{du}\,M_t \log F_{-t}(e^u) = \frac{e^{-\beta t}M_t\,e^u}{F_{-t}(e^u)\,e^{-\beta t}F'_t \circ F_{-t}(e^u)} < e^u
\]
(supposing $u > 0$) remains bounded away from zero. Now $e^{-\beta t}M_t$ converges to something non-zero by the Kesten-Stigum theorem, so $F_t(e^u)\to \infty$ for $u > 0$ by convexity of $F_t$. But then $F_t\circ F_{-t}(e^u) = e^u$ implies $F_{-t}(e^u)\to 1$. All that is left show is that $e^{-\beta t}F'_t \circ F_{-t}(e^u)$ remains bounded. Write $s$ instead of $e^u$. We have \cite{an}
\[
F_t (s) = \big(1 - G(t)\big) s + \int_0^t f \circ F_{t - u}(s)\,dG(u)\:,
\]
so
\begin{eqnarray}
\lefteqn{X_t(s) := e^{-\beta t} \, F'_t(s)}\label{xt}\\
&=&e^{-\beta t}\big(1 - G(t)\big) + \int_0^t f' \circ F_u(s)\,e^{-\beta u}F'_u(s)\,e^{-\beta (t - u)} dG(t - u)\nonumber\\
&=&e^{-\beta t}\big(1 - G(t)\big) + \int_0^t \left(\frac{f' \circ F_u(s)}{\mu} - 1\right) \,X_u(s)\,dG_{\beta}(t - u)\nonumber\\
&& + \:\int_0^t X_u(s)\,dG_{\beta}(t - u)\:,\nonumber
\end{eqnarray}
where we have introduced the measure
\[
G_\beta (t) = \mu \int_0^t e^{-\beta u}\,dG(u)\:.
\]
We now use the final expression in \reff{xt} to obtain
\begin{eqnarray}
\lefteqn{\int_0^t X_u (s)\,dG_{\beta}(t - u) - \int_0^t e^{-\beta u}\,\big(1 - G(u)\big)\,dG_{\beta}(t - u)}\nonumber\\
&=&\int_{u = 0}^t \int_{v = 0}^u \frac{f' \circ F_v(s)}{\mu} \,X_v(s)\,dG_{\beta}(u - v)\,dG_{\beta}(t - u)\nonumber\\
&=&\int_{v = 0}^t \frac{f' \circ F_v(s)}{\mu} \,X_v(s) \int_{u = v}^t dG_{\beta}(u - v)\,dG_{\beta}(t - u) \nonumber\\
&=& \int_{v = 0}^t \frac{f' \circ F_v(s)}{\mu}\,X_v(s)\,dG_{\beta}^{\ast 2}(t - v)\nonumber\:,
\end{eqnarray}
and find by induction and Fubini's theorem,
\begin{eqnarray}
X_t(s) &=& e^{-\beta t}\big(1 - G(t)\big) + \int_0^t e^{-\beta u}\,\big(1 - G(u)\big)\,dU_{\beta}(t - u)\nonumber\\
&&+\: \int_0^t \left(\frac{f' \circ F_u(s)}{\mu} - 1\right) \,X_u(s) \, dU_{\beta}(t - u)\label{diffxt}\:,
\end{eqnarray}
or
\begin{equation}
X_t(s) - X_t(1) = \int_0^t \left(\frac{f' \circ F_{t - u}(s)}{\mu} - 1\right)\,X_{t-u}(s)\,dU_{\beta}(u)\:,\label{zw}
\end{equation}
where
\[
U_{\beta}(t) = \sum_{i = 1}^{\infty} G_{\beta}^{\ast i}(t)
\]
is essentially the renewal measure for $G_{\beta}$. Suppose $X_t \circ F_{-t}(s)$ does not remain bounded. Then we can find a sequence of values $t_1, t_2, \ldots$ tending to infinity such that
\[
1 \leq \lim_{k\to\infty}\int_0^{t_k} \left(\frac{f' \circ F_{t_k-u}\circ F_{-t_k}(s)}{\mu} - 1\right)\,dU_\beta(u)\:.
\]
But $Z_{t-u}/Z_t \to e^{-\beta u}$ almost surely \cite{sch}, hence
\[
\lim_{k\to\infty}F_{t_k-u}\circ F_{-t_k}(s) = \lim_{k\to\infty}\bE\left( F_{-t_k}(s)^{Z_{t_k} Z_{t_k - u}/Z_{t_k}}\right) \leq s^{e^{-\beta u}}\:,
\]
by Jensen's inequality. Therefore,
\[
1 \leq \int_0^{\infty} \left(\frac{f'(s^{e^{-\beta u}})}{\mu} - 1\right)\,dU_\beta(u) \simeq \frac{f''(1)}{\mu}\,\log s
\]
if $s < e^r$ (recall that $f$ has exponential moments up to oreder $r$), because $dU_{\beta}(u)$ is essentially Lebesgue measure plus a term of order $e^{-\alpha u}du$ for some $\alpha \in (0, \beta]$ \cite{st}. But $1 \leq O(\log s)$ is a contradiction for $s$ sufficiently close to 1. The theorem follows.\hfill$\Box$\\

\section{Limit Laws of Luria-Delbr\"uck Processes}
The Luria-Delbr\"uck (LD) distribution arises as the distribution of types in a two-type Bellman-Harris process or, in a narrower sense, as a limiting distribution of types in such a process. If the life-time distribution is exponential and branching is binary, the theory of the LD distribution is essentially complete; see \cite{wa1}. The same can be said if the life-time distribution of cells is exponential, but cells always produce a fixed number $\kappa$ of mutant or non-mutant progeny \cite{asw, ls}. We shall refer to such a process as a $(1 - \rho, \kappa)$-Luria-Delbr\"uck process, where $\rho$ is the probability that upon division, a non-mutant cell produces one non-mutant and $\kappa$ mutant daughter cells. Mutants only produce $\kappa$ mutant progeny. The following theorem has been proved by Leona Schild in her diploma thesis:
\begin{thls} Denote by $L_n$ the number of non-mutants in a $(1 - \rho, \kappa)$-Luria-Delbr\"uck process when population size has reached $n\kappa + 1$, and the process has been started from a single non-mutant individual. Then
\begin{equation}
\frac{L_n}{n^{1 - \rho}} \to L
\end{equation}
almost surely, and
\begin{equation}
L \stackrel{D}{=} B^{1 - \rho} \cdot L_{\vk}\:,\label{DD}
\end{equation}
where $\vk = \kappa^{-1}$, $B$ is $\big((1 - \rho)\vk, \rho\vk\big)$-Beta distributed, and $L_{\vk}$ is $\vk$-biased $(1 - \rho)$-Mittag-Leffler.
\end{thls}
We explain our terms: By a $(1 - \rho)$-Mittag-Leffler distribution we mean a distribution whose Laplace transform is the Mittag-Leffler function
\begin{equation}
E_{1 - \rho}(u) = \sum_{k = 0}^{\infty} \frac{(-u)^k}{\Gamma\big((1 - \rho)k + 1\big)}\:.
\end{equation}
Next, a $\vk$-biasing of a random variable $X$ with Laplace transform $\bE(e^{-uX})$ is a random variable $X_{\vk}$ whose Laplace transform is
\begin{equation}
\frac{\bE(X^{\vk} e^{-uX})}{\bE(X^{\vk})}\:.\label{sb}
\end{equation}
We refer to \cite{ls} for further background. Our goal here is to deduce
\begin{thasw} The Laplace transform of $L$ is
\begin{equation}
\bE(e^{-u L}) = \Gamma(\vk) \sum_{i=0}^{\infty} {-\vk \choose i} \frac{s^i}{\displaystyle \Gamma\big(i(1-\rho)+\vk\big)}\:.\label{lap}
\end{equation}
\end{thasw}
(which has been proved in \cite{asw} in a rather indirect manner) directly from Theorem 7. We need the following
\begin{lemma} \label{ml} Suppose $Y$ is $(1 - \rho)$-Mittag-Leffler distributed. Then its Mellin transform is
\[
\cM Y(s) = \frac{\Gamma(s)}{\Gamma\big((1 - \rho)(s - 1) + 1\big)}\:.
\]
\end{lemma}
{\it Proof}. First observe that
\[
\frac{\Gamma\big((1 - \rho)(k + 1) + 1\big)}{\Gamma\big((1 - \rho) k + 1\big)} \sim k^{1 - \rho}\:,
\]
which implies that $E_{1 - \rho}$ is analytic on all of $\bC$. Also,
\[
\frac{1}{\Gamma(z)} = \frac{1}{2\pi i} \int_{\cH} \frac{e^{\zeta}}{\zeta^z} \,d\zeta \:,
\]
where $\cH$ is a Hankel contour encircling the negative axis in counterclockwise direction. We fix $\ve > 0$, and choose for $\cH = \cH(\ve)$ a vaguely lollipop-shaped figure as follows: It runs along (in fact, just below) the negative axis from $-\infty$ to $-\ve$, runs around the origin on a circle of radius $\ve$ in counterclockwise direction, and then returns to $-\infty$ just above the negative axis. Then we have
\begin{eqnarray}
\lefteqn{\int_{\cH} \frac{e^{\zeta}}{\zeta^{(1 - \rho)k + 1}} \,dt = \int_{- \infty}^{-\ve} \frac{e^{r e^{- i \pi}}}{r^{(1 - \rho)k + 1} e^{- i \pi (1 - \rho)k - i \pi}}\,dr}\nonumber\\
&& + \int_{- \pi}^{\pi} \frac{e^{\ve e^{i \theta}}}{\ve^{(1 - \rho)k + 1} e^{i \theta (1 - \rho)k + i \theta}} \,d\theta + \int_{-\ve}^{- \infty} \frac{e^{r e^{i \pi}}}{r^{(1 - \rho)k + 1} e^{i \pi (1 - \rho)k + i \pi}}\,dr\:,\nonumber
\end{eqnarray}
and the three integrands are $\sim \ve^{- (1 - \rho)k}$ in order of magnitude. It follows that if $\vert u \vert < \ve^{1 - \rho}$, the sequence of functions
\[
\sum_{k = 0}^n (-u)^k \, \frac{e^{r e^{- i \pi}}}{r^{(1 - \rho)k + 1} e^{- i \pi (1 - \rho)k - i \pi}}
\]
is uniformly convergent in $r$, and that of
\[
\sum_{k = 0}^n (-u)^k \,\frac{e^{\ve e^{i \theta}}}{\ve^{(1 - \rho)k + 1} e^{i \theta (1 - \rho)k + i \theta}}
\]
is uniformly convergent in $\theta$. It is therefore permissible to interchange the order of summation and integration in the following chain of equations, and we find
\begin{eqnarray}
\lefteqn{E_{1 - \rho}(u) = \frac{1}{2\pi i} \sum_{k = 0}^{\infty} (-u)^k \int_{\cH} \frac{e^{\zeta}}{\zeta^{(1 - \rho)k + 1}} \,d\zeta}\nonumber\\
&=& \frac{1}{2\pi i} \int_{\cH} \frac{e^{\zeta}}{\zeta} \sum_{k = 0}^{\infty} \frac{(-u)^k }{\zeta^{(1 - \rho)k}} \,d\zeta = \frac{1}{2\pi i} \int_{\cH} \frac{e^{\zeta}}{\zeta + u \zeta^{\rho}} \,d\zeta\:.\label{mlint}
\end{eqnarray}
Now the right-hand side of Equation \reff{mlint} defines an analytic function in $u$ except possibly for those $u$ which belong to the zero set $\cH_0(\ve) := \{u \in \bC: \zeta + u \zeta^{\rho} = 0\ \mbox{for some}\ \zeta \in \cH(\ve)\}$ of the denominator. But one readily checks that $\cH_0(\ve)$ is a Hankel contour around the {\it positive} axis that winds around the origin along a segment of a circle of radius $\ve^{1 - \rho}$, so that if we choose some $\overline{\ve} > \ve$, we have
\[
\int_{\cH(\ve)} \frac{e^{\zeta}}{\zeta + u \zeta^{\rho}} \,d\zeta = \int_{\cH(\overline{\ve})} \frac{e^{\zeta}}{\zeta + u \zeta^{\rho}} \,d\zeta
\]
on a (connected) segment of the open torus $\{u \in \bC: \ve^{1 - \rho} < \vert u \vert < \overline{\ve}^{1 - \rho}\}$. By analytic continuation, then, \reff{mlint} holds for all $u \in \bC \setminus \cH_0(\ve)$, and in particular holds for all $u$ for which $\Re(u) < - \ve^{1 - \rho}$. By Theorem \ref{main}, therefore, we can write
\begin{equation}
\cM Y(s) = - \frac{\Gamma(s)}{4\pi^2} \int_{\beta - i \infty}^{\beta + i \infty}\int_{\cH} \frac{e^{\zeta}}{\zeta + u \zeta^{\rho}} \,d\zeta \, (-u)^{-s}\,du\label{sing}
\end{equation}
if $\beta < - \ve^{1 - \rho}$. Now one may interchange the order of integration (if necessary again after convoluting with $\sqrt{\pi \epsilon} \, e^{-x^2/\epsilon}$) and then apply the Residue theorem to find
\begin{eqnarray}
\lefteqn{\cM Y(s) = - \frac{\Gamma(s)}{4\pi^2} \int_{\cH}\int_{\beta - i \infty}^{\beta + i \infty} \frac{e^{\zeta}}{\zeta + u \zeta^{\rho}}\,(-u)^{-s}\,du\,d\zeta}\nonumber\\
&=& \Gamma(s)\,\frac{1}{2\pi i} \int_{\cH} \frac{e^{\zeta}}{\zeta^{(1 - \rho)s + \rho}}\,d\zeta = \frac{\Gamma(s)}{\Gamma\big((1 - \rho)(s - 1) + 1\big)}\:,\nonumber
\end{eqnarray}
which again follows from Hankel's representation of the reciprocal of the Gamma function. The proof of the lemma is complete.\hfill $\Box$\\

Next we need
\begin{lemma} \label{beta} Suppose $B$ is $\big((1 - \rho)\vk, \rho\vk\big)$-Beta distributed. Then the Mellin Transform of $B^{1 - \rho}$ is
\[
\cM B^{1 - \rho}(s) = \frac{\Gamma(\vk)}{\Gamma\big((1 - \rho)\vk\big)} \frac{\Gamma\big((1-\rho)(s - 1 + \vk)\big)}{\Gamma\big((1-\rho)(s-1) + \vk\big)}\:.
\]
\end{lemma}
{\it Proof}. We have
\[
\sP(B \leq x) = \frac{\Gamma(\vk)}{\Gamma\big((1 - \rho)\vk\big) \Gamma\big(\rho \vk \big)}\int_0^{x} y^{(1 - \rho)\vk - 1} (1 - y)^{\rho\vk - 1}\,dy\:.
\]
Then $\sP(B^{1 - \rho} \leq x) = \sP(B \leq x^{(1 - \rho)^{-1}})$, and
\begin{eqnarray}
\lefteqn{\int_0^1 \frac{d}{dx}\int_0^{x^{(1 - \rho)^{-1}}} y^{(1 - \rho)\vk - 1} (1 - y)^{\rho\vk - 1}\,dy \, x^{s-1}\,dx}\nonumber\\
& = &\int_0^1 \frac{d}{du}\int_0^{u} y^{(1 - \rho)\vk - 1} (1 - y)^{\rho\vk - 1}\,dy \, u^{(1 - \rho)(s-1)}\,du\nonumber\\
& = &\int_0^1 u^{(1 - \rho)\vk - 1} (1 - u)^{\rho\vk - 1}\, u^{(1 - \rho)(s-1)}\,du = \frac{\Gamma(\rho\vk)\,\Gamma\big((1-\rho)(s - 1 + \vk)\big)}{\Gamma\big((1-\rho)(s-1) + \vk\big)}\nonumber\:,\nonumber
\end{eqnarray}
which already finishes the proof. \hfill $\Box$\\

It follows from Lemma \ref{ml} that the Mellin transform of a $\vk$-biased Mittag-Leffler distribution is
\[
\cM Y_{\vk}(s) = \frac{\Gamma(s + \vk)}{\Gamma\big((1 - \rho)(s + \vk - 1) + 1\big)} \frac{\Gamma(1 + (1 - \rho)\vk)}{\Gamma(1 + \vk)}\:.
\]
Together with Lemma \ref{beta} and Theorem \ref{cv}, this implies that the Mellin transform of $L$ satisfies
\[
\cM L(s) = \frac{\Gamma(s + \vk - 1)}{\Gamma\big((1-\rho)(s-1) + \vk\big)}\:.
\]
Specialising now to $s = i + 1$ for arbitrary integer $i \geq 0$, one has for the Laplace transform of $L$,
\begin{eqnarray}
\lefteqn{\frac{\bE(e^{-uL})}{\Gamma(\vk)} = \sum_{i = 0}^{\infty} \frac{\Gamma(i + \vk)}{\Gamma(\vk)\,\Gamma\big((1-\rho)i + \vk\big)}\,\frac{(-u)^i}{i!}}\nonumber\\
&=& \sum_{i = 0}^{\infty} \prod_{j = 1}^{i} \frac{\vk + j - 1}{j}\,\frac{(-u)^i}{\Gamma\big((1-\rho)i + \vk\big)} = \sum_{i=0}^{\infty} {-\vk \choose i} \frac{u^i}{\Gamma\big(i(1-\rho)+\vk\big)}\:.\nonumber
\end{eqnarray}
This was to be proved.\hfill$\Box$\\

{\bf Acknowledgment.} The author thanks Leona Schild for a huge number of stimulating discussions. He is happy to acknowledge financial support from DFG and NWO as part of a Dutch-German research project on random spatial models from physics and biology, and honoured to acknowledge the personal support of Anton Wakolbinger.

\end{document}